\documentclass[11pt]{article}
\usepackage{amssymb}
\usepackage[dvips]{graphics}
\usepackage{theorem}

\newtheorem{theorem}{{\sc Theorem}}[section]

\newtheorem{proposition}[theorem]{{\sc Proposition}}
\newtheorem{definition}{{\sc Definition}}[section]

\newtheorem{corollary}[theorem]{{\sc Corollary}}
\newtheorem{remark}[theorem]{{\sc Remark}}
\newtheorem{qquestion}[theorem]{{\sc Question}}
\newtheorem{pprobleme}[theorem]{{\sc Problem}}
\newtheorem{notation}{{\sc Notation}}

\newenvironment{proof}{Proof}{\nolinebreak $\Box $}

\newenvironment{resume}{\small \begin{center} {\bf R\'esum\'e}
\end{center} \hspace{1cm} \begin{minipage}[t]{13cm} \hspace{.5cm} } 
{\end{minipage} \normalsize }

\begin{document}

\title{Diffeomorphic moment-angle manifolds with different Betti numbers}
\author{Fr\'ed\'eric BOSIO}
\maketitle

\begin{resume}
  Nous construisons ici deux polytopes simples ayant des nombres de Betti 
diff\'erents et auxquels les vari\'et\'es moment-angle associ\'ees sont 
diff\'eomorphes.
\end{resume}

\begin{abstract}
  We find here two simple polytopes whose moment-angle manifolds are 
diffeomorphic and that have different Betti numbers.
\end{abstract}

\section*{Introduction}

  Moment-angle complexes and manifolds are nowadays become a usual object of 
toric topology. They have been introduced by Davis and Januszkiewicz \cite{DJ} 
as "universal spaces" for quasitoric manifolds, and appeared since then in 
several contexts such as toric topology, polyhedral products, intersections of 
quadrics... For an overview of the current theory on moment-angle manifolds and 
related topics, we can refer to~\cite{BP2}.

  The topology of a moment-angle manifold is completely encoded by the 
combinatorial structure of the underlying polytope and we know how to compute 
some topological invariants of a moment-angle manifold, such as cohomology, 
from combinatoric datas on the polytope.

   However, the precise differential structure of a "generic" moment-angle 
manifold is still poorly understood; only several cases have been settled 
(\cite{LdMV}, \cite{BP1}, \cite{BM}, \cite{LdMG}, \cite{CFW}).

  Another challenging problem is to understand when two polytopes give rise to 
"the same" moment-angle manifolds.

  A natural notion of similarity, given differential manifolds, is the notion 
of diffeomorphism. Two polytopes are said {\em{diff-equivalent}} if their 
moment-angle manifolds are diffeomorphic. Despite the apparent strengh of this 
request, a classification of all diff-equivalent polytopes is far out of 
reach, so we shall focus on less ambitious problems. We rather ask which 
combinatorial invariants are kept between diff-equivalent polytopes. 

  The simplest numerical invariant of a polytope is its dimension, and the 
problem of equality of dimensions of diff-equivalent polytopes seems still 
open. Let's just remark that preservation of dimension is equivalent to 
preservation of the number of facets, as the dimension of the moment-angle 
manifold defined by a polytope is the sum of the dimension and the number of 
facets of this one.

  Other important invariants of a polytope are its Betti numbers 
(see~\cite{BP1}).

  Each usual Betti number of a moment-angle manifold, which is well know to be 
a topological invariant, is the sum of some precise Betti numbers of the 
underlying polytope. Hence, if two polytopes have the same Betti numbers, their 
associated moment-angle manifolds have the same usual Betti numbers, and we can 
wonder whether the converse is also true, i.e. if the Betti numbers of the 
moment-angle manifold determine the Betti numbers of the polytope.

  In this paper, we answer negatively to this question, by exhibiting a 
counterexample, i.e. two diff-equivalent polytopes with different Betti numbers.

  To achieve this, we consider multiwedges over neighbourly dual polytopes. For 
such polytopes, the associated moment-angle has a known differential structure, 
which is completely determined by its usual Betti numbers. Indeed, by the work 
of L\'opez de Medrano and Gitler \cite{LdMG}, any such polytope induces a 
connected sums of sphere products as moment-angle manifold, so its differential 
structure is given by the dimensions of the appearing spheres.

  We more precisly consider multiwedges over neighbourly dual polytopes with 
four more facets than their dimension. Indeed, considering multiwedges over 
polytopes with only three more facets than their dimension cannot produce 
counterexample, as in this case, Betti numbers of the moment-angle manifold 
correspond to Betti numbers of the polytope. In the other sense, considering 
multiwedges over polytopes with many more facets than their dimension implies a 
"high" total Betti rank and consequently more numbers to manage. In addition, 
the Gale diagrams of duals ofsuch polytopes take place in high dimensional 
spaces, which increases difficulty. This explains the choice on the kind of 
polytopes we deal with.

  The smallest example of such a polytope is the hexagon. I originally hoped to 
decide whether two diff-equivalent multiwedges over the hexagon must have the 
same Betti numbers. I investigated the cases that, in my opinion, would most 
plausibly provide a counterexample and find none, despite some cases "very 
close to counterexamples". I cannot exclude to have overlooked something. 
Perhaps an exhaustive investigation of all cases is yet possible, which would 
definitely classify the diff-equivalent multiwedges over the hexagon, but this 
is not clear to me due to their abundance.

  In higher dimension, the combinatorics of neighbourly dual polytopes with 
four more facets than their dimension is not unique. There is a broad choice 
of polytopes and we can ask if performing multiwedges with the same multiindex 
over different polytopes can produce counterexamples. This is indeed the case. 
We will consider polytopes for which the difference between Betti numbers can 
be controlled. With suitable choices of a pair of neighbourly polytopes and of 
a multiindex, we will get what we are looking for.

  The counterexample we produce involves $47$-dimensional polytopes with $51$ 
facets, so their common moment-angle manifold has dimension $98$. We can 
suspect this is far from optimal.

  Naturally, one counterexample generates many others, for example by simply 
taking products, we can get as many diff-equivalent polytopes as we want, no 
two of them having the same Betti numbers.

 Another noticeable fact is that our two diff-equivalent polytopes don't have 
equal number of vertices, which proves that the number of vertices is not an 
invariant of diff-equivalence (in contrast, this number, and even the full 
$f$-vector of a polytope, can be recovered from the Betti numbers of this 
polytope ~\cite{BP1} as brought to the author's mind by A. Bahri and T. Panov).

\subsection*{Acknowledgements}

  The author thanks the organisers of the conference "Topology of torus actions 
and applications to geometry and combinatorics", Daejeon, august 2014 
and UMR CNRS 7348 for its financial support.

  The author is grateful to Ivan Limonchenko for pointing out a misprint in the 
original manuscript.

\section{Recalls}

\subsection{Betti numbers of polytopes and cohomology of moment-angle manifolds}

  Given a simple polytope $P$, with facets $F_1 ,..., F_n $, its moment-angle 
manifold can be defined by the quotient operation $P \times T^n / \sim $ where 
$(p,(z_1 ,..., z_n ))$ is identified with $(p,(z'_1 ,..., z'_n ))$ if 
$z_i = z'_i $ whenever $p$ is not on $F_i $ (in other terms, each coordinate of 
$T^n $ generates a rotation that fixes the preimage of the corresponding facet 
of $P$).

  Recall we have the following usual decomposition of the homology of a 
moment-angle manifold, in terms of subsets of facets \cite{B}:

$$H_k (Z_P , {\mathbb Z }) \simeq 
\bigoplus_{{\cal X } \subset {\cal F}} 
\tilde{H}_{k - |{\cal X }| - 1} (P_{\cal X } , {\mathbb Z })$$
A homology class of $Z_P $ is called {\em{induced}} by a subset of facets when 
it is in the image of the reduced homology of this subset.

  We also remark a bigraduation in the homology of the moment-angle manifold, 
by $k$ and $|{\cal X }|$. Fix two integers $p$ and $q$. We note here the Betti 
numbers of $P$ in the following way:
$$b^{p,q} = \bigoplus_{{\cal X } \subset {\cal F} \\ |{\cal X }| = q} 
\dim \tilde{H}_p (P_{\cal X } , {\mathbb Z })$$
This number was standardly noted $b^{p-q+1,2q}$ but the used notation is more 
clear in our context.

\subsection{Balanced configurations}

\begin{definition}
  Let $V$ a $d$-dimensional vector space over ${\mathbb R }$. 
We will call {\em{balanced configuration of points}} in $V$ a configuration 
(i.e. a finite set) of points of $V$ such that:

i) Any $d$ points of the configuration form a basis of $V$.

ii) Given $d-1$ points of the configuration, there are as many points on each 
side of the hyperplane they span.
\end{definition}

  We will note $n$ the number of points of the configuration.

  Remark that with this definition (quite restrictive), the number of points 
cannot have the same parity as $d$.

  For example, the Gale diagram of an even-dimensional neighbourly polytope is 
a balanced configuration of points.

  Here, we're only interested in the combinatorics of such a configuration, 
i.e. the subsets of points whose convex hull contains the origin. Two 
configurations will be called {\em{equivalent}} if they have the same 
combinatorics.

\paragraph{In dimension $2$}

  Let's recall several facts about such configurations:

\begin{proposition}
  In ${\mathbb R }^2 $, a regular polygon centered at the origin with an odd 
number $n$ of sides is a balanced configuration of points. Conversely, any 
balanced configuration of points in a two-dimensional vector space is 
equivalent to such a one.
\end{proposition}

  So, a balanced configuration of points induces a natural cyclic order on its 
points, given by the cyclic order of the polygon. It also induces a distance 
between the points, the distance $d(A,B)$ between two points $A$ and $B$ being 
both their distance on the boundary of the forementioned polygon (considered as 
a graph) and the number of points $C$ such that $0$ is in the convex hull of 
$\{ A , B , C \} $.

\begin{remark}
  Given three points $A,B,C$, $0$ is in their convex hull if and only if 
$d(A,B) + d(B,C) + d(C,A) = n$.
\end{remark}

  Let's give some definitions:

\begin{definition}
  Two points $A$ and $B$ of a balanced configuration of points are said 
{\em{adjacent}} if $d(A,B) = 1$.

  They are said {\em{distant}} if their distance is maximum, i.e. 
$\frac{n-1}{2}$.
\end{definition}

\begin{definition}
  Given two points $A$ and $B$, we note $\phi (A,B)$ the set of points $C$ such 
that $0$ is in the convex hull of $\{ A,B,C \} $.
\end{definition}

  As recalled thereabove $\phi (A,B)$ has $d(A,B)$ elements.

\begin{definition}
  Given a point $A$, the two sets of points on each side of the the line 
spanned by $A$ are called the $A$-classes.

  Both have $\frac{n-1}{2}$ points.
\end{definition}

  Remark that two distant points different from $A$ are one in each $A$-class.

\paragraph{In dimension $3$}

\begin{remark}
  Given a balanced configuration of points in a vector space $V$ 
($\dim V \geq 2$) and a point $x$ of the configuration, the images of the other 
points on the quotient vector space $V / {\mathbb R }x$ yield another balanced 
configuration of points.
\end{remark}

  Hence, given a balanced configuration of points in a $3$-dimensional vector 
space and one of its points $x$, the quotient by ${\mathbb R }x$ induces a 
balanced configuration of points in a $2$-dimensional vector space. The induced 
cyclic order is then called the $x$-cyclic order.

  Let's give some properties of these orders:

\begin{proposition}
  Let $A$ and $B$ two points. Then the $B$-classes for the $A$-order also are 
the $A$-classes for the $B$-order.
\end{proposition}

  Both are the sets of points on each side of the plane spanned by $A$ and $B$, 
hence the result by symmetry.

  The following proposition and its corollary can be useful in determining if 
some convex hulls contain the origin.

\begin{proposition}
  Consider a balanced configuration of points in ${\mathbb R }^3 $, $x$ one of 
its points, and $u$, $v$, $w$ three other points. Then:

i) If $0$ is not in the convex hull of the projections of $u$, $v$ and $w$ on 
${\mathbb R }^3 / {\mathbb R } x$, then $0$ is not in the convex hull of 
$\{ u,v,w,x \} $ in ${\mathbb R }^3 $.

ii)  If $0$ is in the convex hull of the projections of $u$, $v$ and $w$ on 
${\mathbb R }^3 / {\mathbb R } x$, then $0$ is in the convex hull of exactly 
one of $\{ u,v,w,x \} $ and $\{ u,v,w,-x \} $.
\end{proposition}

\begin{proof}
  The point i) is obvious.

  If $0$ is in the convex hull of the projections of $u$, $v$ and $w$ on 
${\mathbb R }^3 / {\mathbb R } x$, then some point $y$ of ${\mathbb R } x$, 
namely $y = \lambda x$, is a positive combination of $u$, $v$ and $w$, and 
$y \neq 0$ by affine independance. Hence, if $\lambda > 0$, then 
$0 = y + \lambda (-x)$ is in the convex hull of $\{ u,v,w,-x \} $ and if 
$\lambda < 0$, then $0 = y + (-\lambda )x$ is in the convex hull of 
$\{ u,v,w,x \} $. And $0$ cannot be on both convex hulls because the 
intersection of the two convex hulls lies on the hyperplane containing 
$u$, $v$ and $w$.
\end{proof}

\begin{corollary}
\label{convhull}
  Assume $0$ is in the convex hull of the projections of $u$, $v$ and $w$ on 
${\mathbb R }^3 / {\mathbb R } x$ and that, for the $u$-cyclic order, $x$ and 
$v$ are on different $w$-classes.

  Then $0$ is in the convex hull of $\{ u,v,w,x \} $.

  Indeed, $-x$ and $v$ are on the same side of the hyperplane spanned by $u$ 
and $w$. So $0$ cannot be in the convex hull of $\{ u,v,w,-x \} $. The 
assertion then derives from the proposition.
\end{corollary}

\subsection{Lexicographic extensions}

  We present here a particular case of a more general construction, well 
explained in ~\cite{P}.

  Given a balanced configuration of points in a $3$-dimensional vector space, 
the following construction produces other ones with two more points:

  Consider three points $u$, $v$, $w$ of the configuration, three signs (i.e. 
$\pm 1$) $s_u $, $s_v $, $s_w $, and two positive real numbers 
$0 < \omega << \epsilon << 1$. Put now 
$p = s_u \cdot u + \epsilon s_v \cdot v + \epsilon ^2 s_w \cdot w$ and 
$q = -p - \omega u - \omega ^2 v - \omega ^3 w$.

  Then adding $p$ and $q$ to the initial configuration yields another balanced 
configuration (see \cite{P}). We will note it 
$[u ^{s_u } , v^{s_v } , w^{s_w }]$.

\paragraph{Lexicographic extensions and dual neighbourly polytopes}

  As we have recalled, if we consider a(n even-dimensional) neighbourly dual 
polytope $P$, the Gale diagram of its dual is a balanced configuration of 
points.

  We also have the converse, i.e. any balanced configuration of points is the 
Gale diagram of a neighbourly polytope. Indeed, if this configuration has $n$ 
points in dimension $d$, it is the Gale diagram of some polytope $P^* $ of 
dimension $n-d-1$ with $n$ vertices. On each side of a hyperplane spanned by 
$d-1$ points of the configuration, there are $\frac{n-d+1}{2}$ points. Hence on 
any side of any vector hyperplane, there is at least $\frac{n-d+1}{2}$. If we 
remove any $\frac{n-d-1}{2}$ points to the configuration, there is at least one 
remaining point on any side of any vector hyperplane. So the origin is in the 
convex hull of the remaining points, which, by properties of Gale diagrams, 
means that any set of $\frac{n-d-1}{2} = \frac{\dim P^* }{2}$ points of 
$P^* $ determins a simplex of $P^* $, i.e. $P^* $ is actually neighbourly.

  Recall also that points of the Gale diagram of the dual $P^* $ of a polytope 
$P$ correspond to facets of $P$.

\subsection{Biflips}

  Consider a simple $d$-polytope $P$. Recall that a flip from $P$ is a passage 
to another simple $d$-polytope $Q$ such that there is a simple $d+1$-polytope 
$T$ having $P$ and $Q$ as disjoint facets and such that there is exactly one 
vertex $v$ of $T$ that neither lays on $P$ nor on $Q$ (see~\cite{T}). There are 
$d+1$ edges of $T$ containing $v$, $p$ of them having their other extremity on 
$P$ and $q$ on $Q$. The flip from $P$ to $Q$ is then called a $(p,q)$-flip, 
with $p+q = d+1$.

  By symmetry, $P$ is then obtained from $Q$ by the inverse $(q,p)$-flip.

  If $p,q \geq 2$, then $P$ and $Q$ have the same number of facets and their 
facets naturally correspond (a facet of $P$, which is the intersection of $P$ 
with a facet of $T$, corresponds to the intersection with $Q$ of the same 
facet).

  We will use the following notation for a flip: The vertices of $P$ adjacent 
to $v$ are the vertices of a simplicial face $F_P $ of $P$, intersection of a 
set $X$ of facets of $P$, called the containing facets. Call $Y$ the set of 
facets of $P$ that meet $F$ on its boundary, which are called the extremal 
facets. (We can notice that containing and extremal facets are inversed in the 
inverse flip).

  The flip will then be noted $(X)[Y]$. The inverse flip from $Q$ to $P$ is 
then $(Y)[X]$.

  Combinatorially, if we unite $X$ with ($Y$ but one element), their 
intersection is a vertex of $P$, whereas if we unite $Y$ with ($X$ but one 
element), their intersection is a vertex of $Q$. All other intersections of 
facets giving vertices of $P$ or $Q$ are the same.

\begin{remark}
  If $P$ is a $d=2d'$ neighbourly dual polytope, then any $(p,q)$-flip from 
$P$, with $p \geq 2$, is a $(d',d'+1)$-flip.

  Indeed, let $(X)[Y]$ such a flip. Then $[Y]$ corresponds to a set of facets 
that do not intersect in $P$, so has at least $d'+1$ elements, and the 
intersection of all elements of $X$ only meets elements of $Y$, so does not 
meet every other facet. It has then at least $d'$ elements by dual 
neighbourlyness.
\end{remark}

\begin{definition}
  We call {\em{biflip}} a pair of flips of the form $(X)[Y]$ and $(Y')[X]$ from 
a polytope $P$ to another polytope $R$, i.e. such that the containing facets of 
the first flip are the extremal facets of the second one.

  Such a biflip will be noted $(X)[Y,Y']$
\end{definition}

  We can notice that the inverse flips from $R$ to $P$ also form a biflip.

\begin{definition}
  Consider a biflip $(X)[Y,Y']$ from a $(d,d+4)$ neighbourly dual polytope. 
Then it is called a {\em{$1$ fixed facet biflip}}, or {\em{$1$ff-biflip}} if 
$X \cup Y \cup Y' $ contains $d+3$ facets, i.e. all but one.
\end{definition}

  Remark that this notion is symmetric, i.e. that the resulting polytope $Q$ is 
also a $(d,d+4)$ neighbourly dual polytope and that $P$ is obtained from $Q$ by 
a $1$ff-biflip (this biflip is then $(X)[Y',Y]$).

\subsection{Wedges and multiwedges}

  Given a simple polytope $P$ and a facet $F$ of $P$, the wedge $W_F P$ of $P$ 
over $F$ is a simple polytope~\cite{KW}.

  It is possible to perfom iterated wedges over a face and even multiwedges by 
performing a fixed number of wedges over each facet (the order of performance 
of wedges are indifferent). Hence, such a mutiwedge is given by an initial 
polytope and a multiindex, a family of natural numbers indexed by the facets of 
this polytope.

  Given a polytope $P$ and a facet $F$, we can compare homology of $Z_P $ and 
of $Z_{W_F P}$. Indeed, if a subset $\cal X $ not containing $F$ induces a 
homology class in $Z_P $, it induces a class of the same bidegree for 
$Z_{W_F P}$, whereas if it contains $F$, it induces in $Z_{W_F P}$ a class in 
which the bidegree is increased by $(1,1)$ (hence the total degree by $2$).

  And this phenomenon passes to multiwedges. In particular, the total Betti 
rank of the moment-angle manifold is not modified by (multi)wedges.

\section{Wedges and biflips}

  We consider here two neighbourly dual polytopes with four facets more than 
their dimension, and who differ from a $1$ff-biflip. We compute the difference 
between the subsets inducing homology in their moment-angle manifolds. Let's
recall that on the boundary of a $d=2d'$-dimensional neighbourly dual polytope, 
a union of facets (except all or none) has the homotopy type of a wedge of 
spheres of dimension $(d'-1)$.

  Assume we have a $d=2d'$-dimensional neighbourly dual polytope $P$ with $d+4$ 
facets. Let ${\cal F } $ the set of its facets and ${\cal F }'$ any nonempty 
subset of ${\cal F }$ with at most $d'$ elements. As $P$ is neighbourly dual, 
there is a vertex of $P$ on the intersection of all the facets belonging to 
${\cal F }'$. The union of all these facets is starshaped on such a vertex, so 
${\cal F }'$ does noy contribute to the homology of $Z_P $. By duality, neither 
does its complement, so all subsets of ${\cal F }$ inducing homology in $Z_P $ 
(except $\emptyset $ and ${\cal F }$ itself) have $d'+1$, $d'+2$ or $d'+3$ 
elements.
 
  Assume now we have two $(d,d+4)$ neighbourly dual polytopes $P_1 $ and 
$P_2 $, with $d=2d'$, obtained from each other by a $1$ff-biflip, namely 
$(X , Y_1 , Y_2 )$. Call then $Y = Y_1 \cap Y_2 $.

  Now, $X$ has $d'$ elements, $Y_1 $ and $Y_2 $ have $d'+1$ elements and, as 
the flip is $1$ff, $X \cup Y_1 \cup Y_2 $ has $d+3$ elements, so 
$Y_1 \cup Y_2 $ has $d'+3$ and $Y$ has $d'-1$ elements. There are then exactly 
two elements in each that are not in the other. Note $G_1 $, $H_1 $ the facets 
that are in $Y_1 $, not in $Y_2 $ and $G_2 $, $H_2 $ the facets that are in 
$Y_2 $, not in $Y_1 $ and $G$ the facet which is neither in $X$, $Y_1 $ nor 
$Y_2 $.

  We are looking for the subsets of $d'+1$ facets that induce homology in one 
of the two moment-angle manifolds, not in the other. A subset of $d'+1$ facets 
induces homology in $Z_{P_1 } $ or $Z_{P_2 }$ exactly when the intersection of 
all its members is empty in the corresponding polytope. So, if such a subset 
induces homology in $Z_{P_1 }$ and not in $Z_{P_2 }$, the intersection of its 
members must be destroyed by the $1$ff-biflip.

  The first flip $(X)[Y_1 ]$, destroyes the intersection of $X \cup \{ F \} $ 
for any $F$ in $Y_1 $, and no other intersection of sets with $d'+1$ elements. 
The second flip $(Y_2 )[X]$ destroyes the intersection of $Y_2 $ and no other 
intersection of sets with $d'+1$ elements, and it reintroduces intersection of 
$X \cup F$ if $F$, assumed in $Y_1 $, is also in $Y_2 $.

  Hence, there are only three susbets of $d'+1$ facets whose members intersect 
in $P_1 $ and not in $P_2 $, namely $Y_2 $, $X \cup \{ G_1 \} $ and 
$X \cup \{ H_1 \} $. By symmetry, there are three susbets of $d'+1$ facets 
whose members intersect in $P_2 $ and not in $P_1 $, namely $Y_1 $, 
$X \cup \{ G_2 \} $ and $X \cup \{ H_2 \} $.

  Let's now determine the subsets of $d'+2$ facets whose contribution to 
moment-angle homology changes after the biflip. We can remark that the 
contribution to moment-angle homology of such a subset only depends on the 
number of its subsets with $d'+1$ that do not intersect in the polytopes. So 
if such a subset contributes more to the homology of $Z_{P_1 }$ than to the one 
of $Z_{P_2 }$, it must contain one of the three subsets mentioned thereabove. 
By duality, so does its complement, and one of them must contain $Y_2 $, the 
other must contain $X \cup \{ G_1 \} $ or $X \cup \{ H_1 \} $, at least $X$.

  So there are only three pairs of possibilties, namely $Y_2 \cup \{ G_1 \}$ 
and $X \cup \{ H_1 , G \} $, $Y_2 \cup \{ H_1 \}$ and $X \cup \{ G_1 , G \} $, 
$Y_2 \cup \{ G \} $ and $X \cup \{ G_1 , H_1 \} $.

  Note the coherence: $X \cup \{ G_1 , H_1 \} $ actually contains two new 
subsets of $d'+1$ elements that do not inetersect, but as they globally 
intersect in $P_1 $, they only contribute for one dimension in homology of 
$Z_{P_2 }$.

  To sum up, here are the changes in moment-angle homology:

\begin{itemize}

\item
$Y_2 $, $X \cup \{ G_1 \} $ and $X \cup \{ H_1 \} $ contribute for one 
dimension in $H^{d+1}(Z_2 )$, not in $H^{d+1}(Z_1 )$.

\item 
$Y_2 \cup \{ G_1 \}$, $X \cup \{ H_1 , G \} $, $Y_2 \cup \{ H_1 \}$, 
$X \cup \{ G_1 , G \} $, $Y_2 \cup \{ G \}$, $X \cup \{ G_1 , H_1 \} $ 
contribute for one dimension more in $H^{d+2}(Z_2 )$ than in $H^{d+2}(Z_1 )$.

\item
(By duality) $X \cup \{ G_1 , H_1 , G \} $, $Y_2 \cup \{ H_1 , G \} $ and 
$Y_2 \cup \{ G_1 , G \} $ contribute for one dimension in $H^{d+3}(Z_2 )$, not 
in $H^{d+3}(Z_1 )$.

\item
$Y_1 $, $X \cup \{ G_2 \} $ and $X \cup \{ H_2 \} $ contribute for one 
dimension in $H^{d+1}(Z_1 )$, not in $H^{d+1}(Z_2 )$.

\item 
$Y_1 \cup \{ G_2 \}$, $X \cup \{ H_2 , G \} $, $Y_1 \cup \{ H_2 \}$, 
$X \cup \{ G_2 , G \} $, $Y_1 \cup \{ G \}$, $X \cup \{ G_2 , H_2 \} $ 
contribute for one dimension more in $H^{d+2}(Z_1 )$ than in $H^{d+2}(Z_2 )$.

\item
$X \cup \{ G_2 , H_2 , G \} $, $Y_1 \cup \{ H_2 , G \} $ and 
$Y_2 \cup \{ G_2 , G \} $ contribute for one dimension in $H^{d+3}(Z_1 )$, not 
in $H^{d+3}(Z_2 )$.

\end{itemize}

  If we now consider a finite sequence of sufficiently distant $1$ff-biflips, 
so that no two subsets of facets inducing change in homology in diffent biflips 
are close to each other, then the global change after the sequence of biflips 
is the addition of each change.

  When considering multiwedges over such polytopes, with the same multiindex, 
the changes are the same concerning dimensions (as vector spaces), but degrees 
and bidegrees on which the changes occur depend on the multiindex.

  So we ask whether there exists a multiindex for which the bigraded Betti 
numbers of the polytopes are different, but the Betti numbers of moment-angle 
manifolds are the same. We have found such a possibilty with four biflips, with 
a multiindex satisfying the values given in the table, 
section~\ref{sectncounterex}.

\section {Lexicographic extension and biflips}

  We show here how to construct pairs of polytopes with the required property. 

\begin{definition}
  Consider a balanced configuration of points in ${\mathbb R }^3 $, $x_1 $ one 
of these points, and $x_2 $, $x_3 $ two points that are adjacent for the 
$x_1 $-cyclic order.

  Then, the lexicographic extensions $[x_1 ^{(+)} , x_2 ^{(+)} ,  x_3 ^{(+)}]$ 
and $[x_1 ^{(+)} , x_3 ^{(+)} ,  x_2 ^{(-)}]$ are called {\em{cousin 
extensions}}.

  Calling $y$ the common distant point of $x_2 $ and $x_3 $ for the 
$x_1 $-cyclic order, the points $x_1 $, $x_2 $, $x_3 $, $y$ and the two ones we 
add, $p$ and $q$, are called the {\em{special points}} of the extensions. The 
other ones will be called {\em{indifferent}}.
\end{definition}

\underline{Warning:} The two extensions are not exactly symmetric in $x_2 $, 
$x_3 $, due to the sign $-$ in the second one. The first extension will be 
called the left extension, the second the right extension.

  We assume here that we have an even-dimensional $(d,d+4)$ neighbourly dual 
polytope $P$. The Gale diagram of its dual yields a balanced configuration of 
$d+4$ points in ${\mathbb R }^3 $.

  Consider two cousin extensions related to three points $x_1 $, $x_2 $ and 
$x_3 $ (corresponding to facets of $P$), and the associated polytopes 
$P_l $ (left extension) and $P_r $ (right extension). Then:

\begin{proposition}
\label{bonsbifl}

  The polytopes $P_l $ and $P_r $ are obtained from each other by a 
$1$ff-biflip.

\end{proposition}

\begin{proof}
  We just have to analyse accurately the differences between these two 
polytopes.

  In the sequel, we identify points of balanced configurations and facets of 
neighbourly polytopes.

  We begin with a definition:

\begin{definition}
 Let's place under the conditions of the proposition. A set of four points of 
our extensions will be called a {\em{changing quatuor}} if it contains $0$ 
in its convex hull for one extension, not for the other.
\end{definition}

  To compute the changing quatuors, we look at the modifications of the cyclic 
orders when joining $p$ and $q$ to our facets.

\begin{itemize}

\item
  Let $F$ be an indifferent facet, and let's determine the $F$-cyclic order of 
the extensions. We just have to place correctly $p$ anq $q$. We see that $p$ is 
infinitesimally close to $x_1 $, so $q$ infinitesimally close to the opposite 
of $x_1 $, hence between the two distant facets of $x_1$ for the $F$-cyclic 
order. So the positions of $q$ are the same in both extension. Furthermore, $p$ 
is in the $x_1 $-class containing $x_2 $ for the first extension and in the 
$x_1 $-class containing $x_3 $ for the second one. Now, as $F$ is not $y$, 
$x_2 $ and $x_3 $ are in the same $F$-class for the $x_1 $-cyclic order. So 
they also are in the same $x_1 $-classe for the $F$-cyclic order.

  Hence, the $F$-cyclic orders are equal for both extensions.

\includegraphics{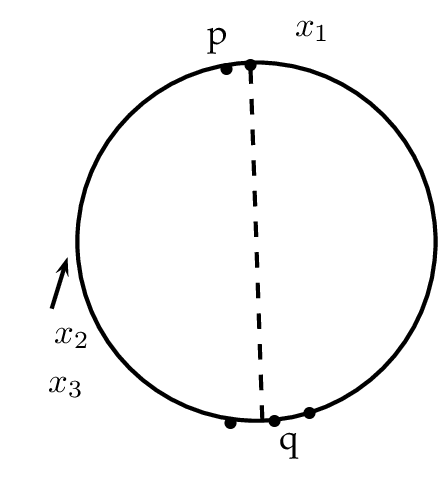}

  Let's now examine the special facets.

\item
  Consider first $y$. Still $q$ is between the distant facets to $x_1 $ for the 
$y$-cyclic order, but $x_2 $ and $x_3 $ are then on opposite $x_1 $-classes for 
the $y$-cyclic order, which implies that, for the two $y$-cyclic orders, $p$ is 
each time adjacent to $x_1 $, but once on each side.

\includegraphics{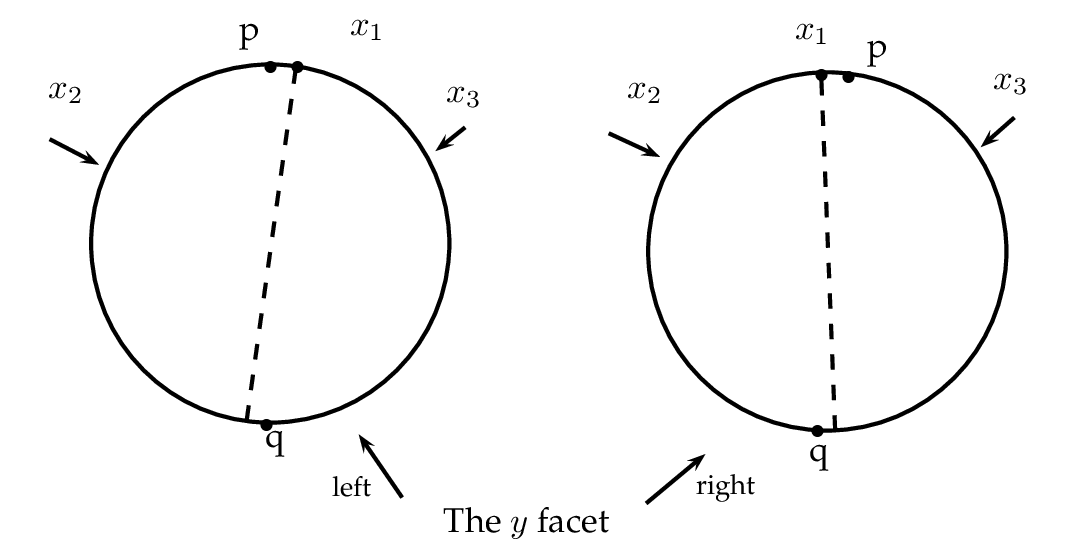}

  Consider $x_2 $. Still $q$ is between the distant facets to $x_1 $ for the 
$x_2 $-cyclic order, and $p$ is adjacent to $x_1 $. We have, for the left 
extension, modulo $x_2 $, $p \equiv x_1 + \epsilon ^2 x_3 $ so $p$ is in 
the $x_1 $-class containing $x_3 $. For the right extension, we have, 
modulo $x_2 $, $p \equiv x_1 + \epsilon x_3 $ so $p$ is also in the 
$x_1 $-class containing $x_3 $.

  Hence, the $x_2 $-cyclic orders are the same for both extensions.

\includegraphics{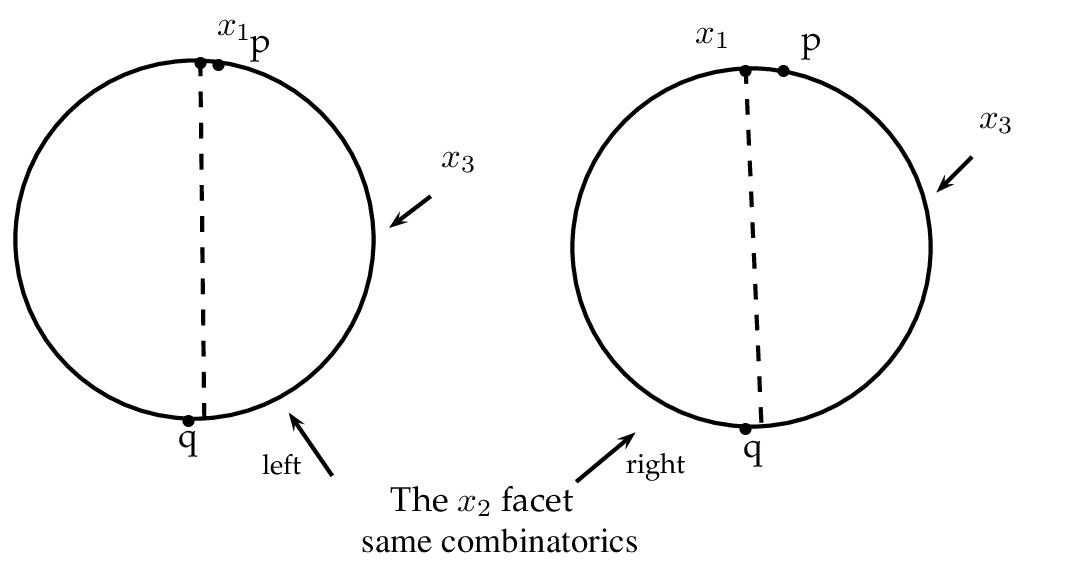}

\item
  Consider $x_3 $. Still $q$ is between the distant facets to $x_1 $ for the 
$x_3 $-cyclic order, and $p$ is adjacent to $x_1 $. We have, for the left 
extension, modulo $x_3 $, $p \equiv x_1 + \epsilon x_2 $ so $p$ is on $x_2 $'s 
side from $x_1 $. For the right extension, we have, modulo $x_3 $, 
$p \equiv x_1 - \epsilon ^2 x_2 $ so $p$ is,from $x_1 $, on the side opposite 
to $x_2 $. Hence, $x_1 $ and $p$ are inversed in the $x_3 $-cyclic orders of 
the extensions.

\includegraphics{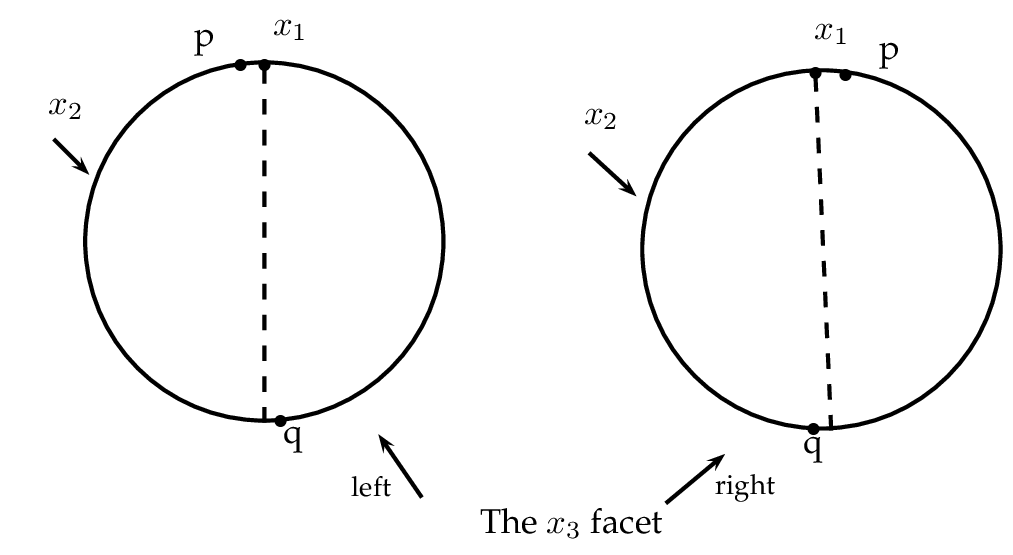}

\item
  Consider $x_1 $. For the left extension, we have, modulo $x_1 $,
$p \equiv \epsilon x_2 + \epsilon ^2 x_3 $, so is adjacent to $x_2 $, in 
the $x_2 $-class containing $x_3 $, i.e. between $x_2 $ and $x_3 $. For the 
right one, $p \equiv \epsilon x_3 - \epsilon ^2 x_2 $, so is adjacent to 
$x_3 $, on the $x_3 $-class not containing to $x_2 $. Whereas $q$ is in both 
cases infinitely close to $-p$, i.e. close to the opposite of $x_2 $ in the 
left extension and to the opposite of $x_3 $ in the right one. For the cyclic 
orders, both corresponding elements are adjacent to $y$, one on each side 
depending on the extension.

\includegraphics{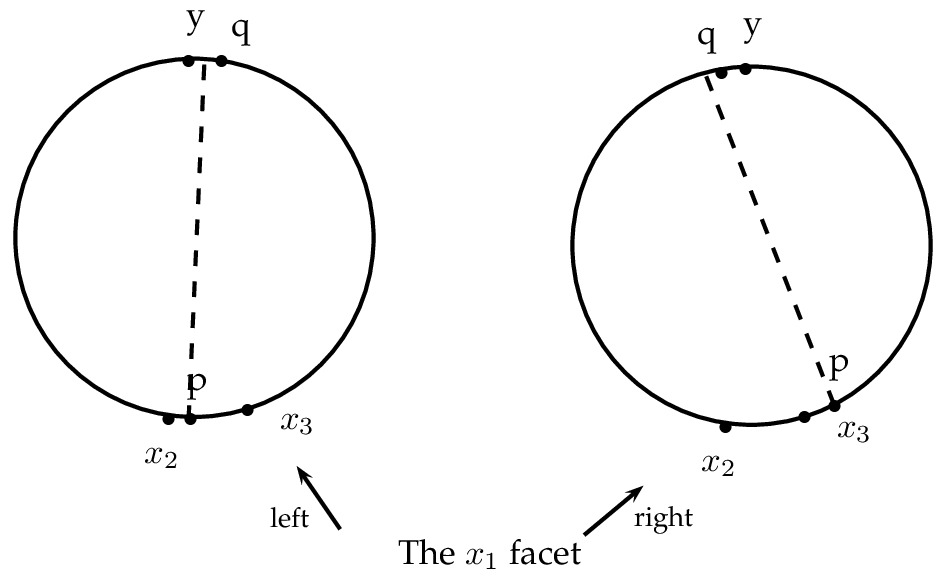}

\item
  Consider $p$. As $p$ is infinitesimally close to $x_1$, the $p$-cyclic order 
on the initial facets but $x_1 $ is the same as the initial $x_1 $-cyclic 
order. We have, modulo $p$, for the left extension, 
$x_1 \equiv x_1 - p = -\epsilon x_2 - \epsilon ^2 x_3 $, so is infinitesimally 
close to $-x_2 $. In the right one, 
$x_1 \equiv x_1 - p = -\epsilon x_3 + \epsilon ^2 x_2 $, so is infinitesimally 
close to $-x_3 $. For the $p$-cyclic order, $x_1 $ is each time adjacent to 
$y$, once on each side.

  We also have, in both cases, modulo $p$, $q \equiv q+p$, which, in the first 
extension, equals 
$-\omega x_1 - \omega ^2 x_2 - \omega ^3 x_3 \equiv 
(\epsilon - \omega) x_2 + (\epsilon ^2 - \omega ^3 ) x_3 $, so $q$ is 
infinitesimally close to $x_2 $, on $x_3 $'s side, i.e. between both, and, in 
the second extension, equals 
$-\omega x_1 - \omega ^2 x_2 - \omega ^3 x_3 \equiv 
(\epsilon - \omega) x_3 - (\epsilon ^2 + \omega ^3 ) x_2 $, so $q$ is 
infinitesimally close to $x_3 $, on the opposite of $x_2 $'s side. So for the 
$p$-cyclic order, $x_1 $ is each time adjacent to $y$, once on each side.

\includegraphics{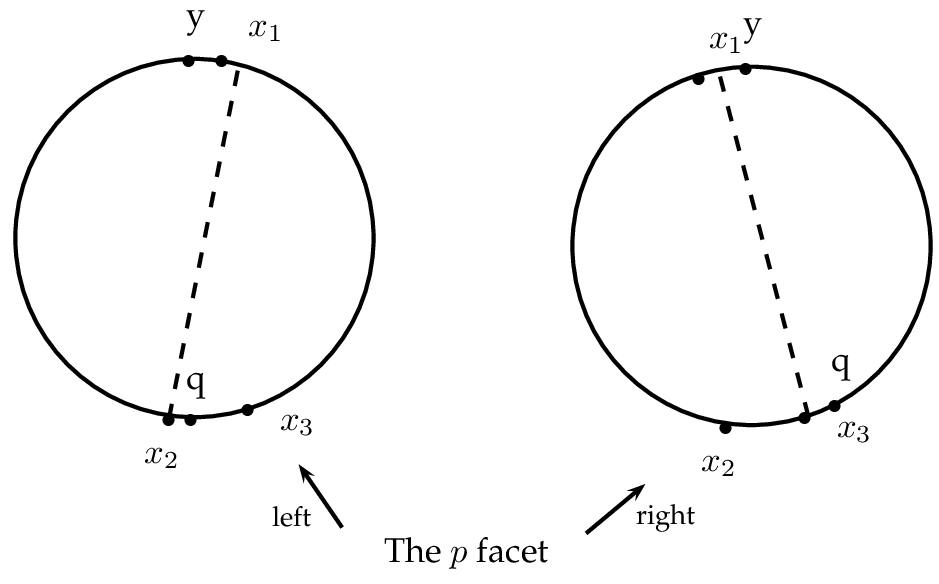}

\item
  Let's finally consider $q$. As $q$ is infinitesimally close to $-x_1 $, the 
cyclic order on the initial facets but $x_1 $ is the same as the initial 
$x_1 $-cyclic order.

We have, in the left extension, 
$q = -(1 + \omega )x_1 - (\epsilon + \omega ^2 ) x_2 - 
(\epsilon ^2 + \omega ^3 ) x_3 $. So, modulo $q$, we have 
$(1+ \omega )x_1 \equiv (1+ \omega )x_1 + q = 
- (\epsilon + \omega ^2 ) x_2 - (\epsilon ^2 + \omega ^3 ) x_3 $. So, for the 
$q$-cyclic order, $1$ is infinitesimally close to $-x_2 $.

  In the right extension, we have 
$q = -(1 + \omega )x_1 - (\epsilon + \omega ^2 ) x_3 + 
(\epsilon ^2 - \omega ^3 ) x_2 $. So, modulo $q$, we have 
$(1+ \omega )x_1 \equiv (1+ \omega )x_1 + q = 
- (\epsilon + \omega ^2 ) x_3 + (\epsilon ^2 - \omega ^3 ) x_3 $. So, for the 
$q$-cyclic order, $1$ is infinitesimally close to $-x_3 $.

  So, for the $q$-cyclic order, $x_1 $ is each time adjacent to $y$, once on 
each side.

  We have, modulo $q$, in the left extension, 
$p \equiv p + \frac{1}{1 + \omega } q = 
(1 - \frac{1 + \omega }{1 + \omega }) x_1 + 
(\epsilon - \frac{\epsilon + \omega ^2 }{1 + \omega }) x_2 +
(\epsilon ^2 - \frac{\epsilon ^2 + \omega ^3 }{1 + \omega }) x_3 = 
\frac{\omega }{1 + \omega } 
\left[ (\epsilon - \omega ) x_2 + (\epsilon ^2 - \omega ^2 ) x_3 \right] = 
\frac{\omega (\epsilon - \omega }{1 + \omega } 
(x_2 + (\epsilon + \omega ) x_3 )$.

  For the $q$-cyclic order, $p$ is between $x_2 $ and $x_3 $.

  In the right one, modulo $q$, we have 
$p \equiv p + \frac{1}{1 + \omega } q = 
(1 - \frac{1 + \omega }{1 + \omega }) x_1 + 
(\epsilon - \frac{\epsilon + \omega ^2 }{1 + \omega }) x_3 -
(\epsilon ^2 - \frac{\epsilon ^2 - \omega ^3 }{1 + \omega }) x_2 =
\frac{\omega }{1 + \omega } 
\left[ (\epsilon - \omega ) x_3 - (\epsilon ^2 + \omega ^2 ) x_2 \right] $.
For the $q$-cyclic order, $p$ is adjacent to $x_3 $, on the opposite side to 
$x_2 $.

  So, for the $q$-cyclic order, $p$ is each time adjacent to $x_3 $, once on 
each side.

\includegraphics{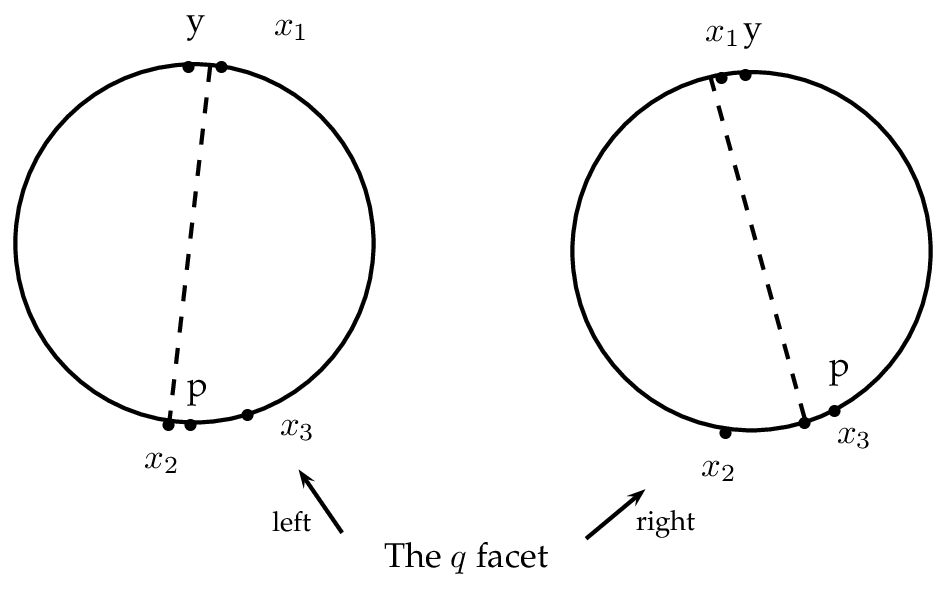}
\end{itemize}

  To sum up, the only inversions in the cyclic orders occuring when switching 
the extensions are:

\begin{itemize}
\item For the $x_1 $-cyclic order, $x_3 $ with $p$ and $y$ with $q$.
\item For the $x_3 $-cyclic order, $x_1 $ with $p$.
\item For the $y$-cyclic order, $x_1 $ with $p$.
\item For the $p$-cyclic order, $x_3 $ with $q$ and $x_1 $ with $y$.
\item For the $q$-cyclic order, $x_3 $ with $p$ and $x_1 $ with $y$.
\end{itemize}

  We now can deduce that the polytopes are obtained from each other by a 
$1$ff-biflip.

  Consider the following sets of facets:

  In the initial $x_1 $ cyclic order, $x_3 $ detrmines two classes, namely 
${\cal X }$, containing $x_2 $ and ${\cal X }'$, containing $y$. Put 
${\cal Y } = {\cal X }' \setminus\{ y \} $, 
${\cal Y }_1 = {\cal Y } \cup \{ p,y \} $ and 
${\cal Y }_2 = ({\cal Y } \cup \{ x_3 , x_1 \} )$. Let's check that $P_2 $ is 
obtained from $P_1 $ by the biflip ${(\cal X })[{\cal Y }_1 , {\cal Y }_2 ]$, 
which is a $1$ff-biflip because $q$ is the only facet neither belonging to 
${\cal X }$, ${\cal Y }_1 $ nor ${\cal Y }_2 $.

  The complement of ${\cal X } \cup {\cal Y }_1 $ contains $x_1 $, $x_3 $ and 
$q$. The complement of ${\cal X } \cup {\cal Y }_2 $ contains $y$, $p$ and $q$.

  Passing to Gale diagrams, we then have to show that sets of the form 
$\{ x_1 , x_3 , q , F \} $ are changing quatuors, $0$ being in their convex 
hull when $F$ is in ${\cal X }$ for the left extension, and in ${\cal Y }_1 $ 
for the right one, that sets of the form $\{ y , p , q , F \} $ are changing 
quatuors, $0$ being in their convex hull when $F$ is in ${\cal Y }_2 $ for the 
left extension, and in $X$ for the right one, and that there are no other 
changing quatuor.

  In the left extension, we have, $\phi _{x_1 } (x_3 , q) = {\cal X }$. We 
deduce that, if $F$ is in ${\cal Y }_1 $, then $0$ is not in the convex hull of 
$\{ x_1 , x_3 , q , F \}$.

  In the right extension, $\phi _{x_1 } (x_3 , q) = {\cal Y }_1 $. So, if $F$ 
is in ${\cal X }$, then $0$ is not in the convex hull of 
$\{ x_1 , x_3 , q , F \}$.

  In the left extension, we have, $\phi _p (y, q) = {\cal Y }_2 $. We 
deduce that, if $F$ is in ${\cal X }$, then $0$ is not in the convex hull of 
$\{ y , p , q , F \}$.

  In the right extension, $\phi _{p} (y,q) = {\cal X }$. So, if $F$ 
is in ${\cal Y }_2 $, then $0$ is not in the convex hull of 
$\{ y , p , q , F \}$.

  Consider the $x_3 $-classes in the $q$-cyclic order for the right extension. 
We have ${\cal X } \cup \{ x_1 \} $ and 
${\cal Y } \cup \{ p,y \} = {\cal Y }_1 $. Hence, by corollary~\ref{convhull}, 
if $F$ is in ${\cal Y }_1 $ then $0$ is actually in the convex hull of 
$\{ x_1 , x_3 , q , F \} $ for the right extension.

  Consider the $x_3 $-classes in the $q$-cyclic order for the left extension. 
We have ${\cal X } \cup \{ p \} $ and 
${\cal Y } \cup \{ x_1 ,y \} = {\cal Y }_1 $. Hence, by 
corollary~\ref{convhull}, if $F$ is in $X$ then $0$ is actually in the convex 
hull of $\{ x_1 , x_3 , q , F \} $ for the left extension.

  Consider the $y$-classes in the $q$-cyclic order for the right extension. 
We have ${\cal X } \cup \{ x_1 \} $ and 
${\cal Y } \cup \{ p,x_3 \} $. Hence, by corollary~\ref{convhull}, if $F$ is in 
$X$ then $0$ is actually in the convex hull of $\{ y , p , q , F \} $ for 
the right extension.

  Consider the $y$-classes in the $q$-cyclic order for the left extension. 
We have ${\cal X } \cup \{ p \} $ and 
${\cal Y } \cup \{ x_1 , x_3 \} = {\cal Y }_2 $. Hence, by 
corollary~\ref{convhull}, if $F$ is in ${\cal Y }_2 $ then $0$ is actually in 
the convex hull of $\{ y , p , q , F \} $ for the left extension.

  We then have seen that the required changes between the two extensions occur. 
Let's now verify there are no other ones.

  Consider four points. If neither $p$ nor $q$ is not one of these points, 
there is nothing to show. If $p$ is one of these points, not $q$, then $0$ is 
in their convex hull if and only if this is also true when replacing $p$ by 
$x_1 $, and this in both extensions. Hence, a set of four points not containing 
$q$ cannot be a changing quatuor.

  So we can consider $q$ and three other points. If neither of them is $x_1 $ 
nor $p$, then $0$ is in their convex hull if and only if this is also true when 
replacing $q$ by $-x_1 $, and this in both extensions. Such a set cannot be a 
changing quatuor.

  If (at least) two of these points are nonspecial or $x_2 $, then they induce 
the same cyclic orders in both extensions and, by corollary~\ref{convhull}, we 
can conclude that this set is not a changing quatuor.

  Hence a changing quatuor must contain $q$ and at least two points among 
$x_1 $, $x_3 $, $p$, $y$. If it contains $x_1$, $x_3 $ and $q$ or $p$, $y$ and 
$q$, we have prouved it is a changing quatuor.

  Only remains the case of $q$, ($x_1 $ or $x_3 $), ($y$ or $p$) and a last 
point.

  Now, $\phi _q (x_1 , y)$ equals $\{ p \} $ for the left extension and 
$\{ x_3 \} $ for the right one. This induces no new changing quatuor.

  Also, $\phi _q (x_3 , p)$ equals $\{ y \} $ for the left extension and 
$\{ x_1 \} $ for the right one. This induces no new changing quatuor.

  Consider $\{ x_3 , y , q , F \}$, $F$ not being $x_1 $ nor $p$. As the only 
inversion in the $x_3 $-cyclic orders are $1$ with $p$, $F$ is in 
$\phi _{x_3 } (y,q)$ for both extensions or for none. Also, the $F$-cyclic 
orders are equal in both extensions, so, by corollary~\ref{convhull}, 
$\{ x_3 , y , q , F \}$ is not a changing quatuor.

  Consider $\{ x_1 , p , q , F \}$, $F$ not being $x_3 $ nor $y$. For the first 
extension, we have $\phi _q (x_1 , p) = X \cup \{ y \} $ abd for the second, 
$\phi _q (x_1 , p) = X \cup \{ x_3 \} $. So $F$ is in 
$\phi _q (x_1 , p)$ for both extensions or for none. As the $F$-cyclic 
orders are equal in both extensions, we again have, by 
corollary~\ref{convhull}, $\{ x_1 , p , q , F \}$ is not a changing quatuor.

  All the cases have been examined. So we are done.
\end{proof}

  We now consider multiple pairs of extensions, i.e. we start with a 
neighbourly dual polytope or a balanced configuration of points and we perform 
several pairs of cousin extensions. We can ask if they can be performed 
together.

\begin{proposition}

  Asuume we have a balanced configuration of points, and a family 
$ce_1 ,..., ce_k $ of cousin extensions of this configuration. Then, if their 
special facets (points) are all different, we can get two extensions by 
performing them in any order and any sense (for each pair or cousins, we can 
choose the left or right one for the first extersion and the other for the 
second). The resulting extensions do not depend on the order in which the 
extensions have been performed (i.e. in some sense, these extensions commute), 
it just depends on the choices done for the senses.

  And the resulting polytopes are obtained by each other by a sequence of $k$ 
consecutive $1$ff-biflips.
\end{proposition}

\begin{proof}
  To verify that the final cyclic orders do not depend on the order of 
performance of the extensions, we just have to see that there is no ambiguity 
on the combinatorics of the polytope, provided the different $\epsilon $'s and 
$\omega $'s are small enough, independantly of their relative sizes.

  Indeed, consider four points for whose we try to determine if $0$ is in their 
convex hull.

\begin{enumerate}

\item
  If the added points among these four ones come from a unique extension, there 
is no ambiguity.

\item
  If among the four points there is some $p$ (resp. $q$) without the 
corresponding $q$ (resp. $p$) nor $x_1 $, then the considered point can be 
replaced by the corresponding $x_1 $ (resp. $-x_1 $) so that "its" extension 
can be peformed at any time.

\item
  If among the four points there is some $p$ or $q$ accompanied with the other 
without the corresponding $x_1 $, then $p$ can be replaced by "its" $x_1 $ and 
$q$ by "its" $x_2 $ or $x_3 $ depending only on the sense of this extension.

\item
 If among the four points there is some $p$ (resp. $q$) accompanied with the 
corresponding $x_1 $ but not with "its" $q$ (resp. $p$), then we can replace 
$p$ (resp. $q$) by $x_2 $ or $x_3 $ (resp. $-x_2 $ or $-x_3 $) depending only 
on the sense of this extension.

\item
  If finally, we have $p$, $q$ and $x_1 $ in the four points, then there is no 
ambiguity because there is only one facet remaining (if it is some $p$ or $q$, 
we are in fact in the case 2).

\end{enumerate}

  Hence, the order of performance of the extensions is indifferent. We now just 
have to check that the two final polytopes are obtained from each other by the 
correct sequence of biflips.

  We can proceed by induction on $k$. Assume this is true for $k$ cousin 
extensions, and consider $k+1$ ones, giving two polytopes $T_{k+1}$ and 
$T'_{k+1}$. Consider then the polytope $T''_{k+1}$ obtained by changing only 
the last choice on the family of extension giving $T'_{k+1}$ (making this 
choice the same as the one giving $T_{k+1}$). Then, by what precedes, 
$T'_{k+1}$ and $T''_{k+1}$ are obtained from each other by a $1$ff-biflip. 
Now, this last extension could also have been performed at first, so that 
$T_{k+1}$ and $T''_{k+1}$ are both obtained from the same polytope by a 
sequence of $k$ consecutive $1$ff-biflips with different special facets. By 
induction assumption, they are obtained from each other by a sequence of $k$ 
consecutive $1$ff-biflips. So $T_{k+1}$ and $T'_{k+1}$ are obtained from each 
other by a sequence of $k+1$ consecutive $1$ff-biflips, which concludes the 
induction.
\end{proof}

\section{Explicit counterexample}
\label{sectncounterex}

  We give here an explicite example of two polytopes with diffeomorphic 
moment-angle manifolds and different Betti numbers.

  Consider the dual cyclic polytope $P = C ^* _{16,20}$, with the natural 
cyclic order on its facets.

  The Gale diagram of its dual is a balanced configuration of points in 
dimension $3$. If we choose the (point corresponding to) facet $1$, it induces 
a cyclic order on the other ones. For this cyclic order, facets $2$ and $4$ are 
adjacent, their common distant facet being $3$.

  We now produce four pairs of counsin extensions, the first one with 
$x_1 = 1$, $x_2 = 2$ and $x_3 = 4$. The special facets are then, in addition 
with the two new ones, $1$, $2$, $3$ (corresponding to $y$) and $4$.

  The three other pairs of extensions are obtained by shifting one, two and 
three times the facets by $5$, hence no two pairs among them have any common 
special facet.

  So we can consider extensions, and associated polytopes, obtained by 
performing these counsin extensions together. We will note these two polytopes 
$P'_1 $ and $P'_2 $

  For  reasons of aesthetics, we alternate the senses of the extensions, i.e., 
for $P'_1 $ we consider the left extensions in the first and third pair and the 
right extension for the second and forth ones, and we inverse for $P'_2 $. So, 
we easily see that shifting by $5$ the facets induces an isomorphism between 
$P'_1 $ and $P'_2 $.

  Now our counterexample $P_1 $ and $P_2 $ will be obtained by multiwedges over 
respectively $P'_1 $ and $P'_2 $, the indices of the multiwedges being the same 
and given by the following table:

$$\begin{array} {|c||c|c|c|c|c|c|c|c|c|c|}
\hline
\mbox{Facet} & 1 & 2 & 3 & 4 & 5 & 6 & 7 & 8 & 9 & 10 \\
\hline
\mbox{index} & 0 & 3 & 1 & 2 & 0 & 1 & 0 & 0 & 2 & 2 \\
\hline \hline
\mbox{Facet} & 11 & 12 & 13 & 14 & 15 & 16 & 17 & 18 & 19 & 20 \\
\hline
\mbox{index} & 1 & 0 & 0 & 1 & 0 & 0 & 3 & 1 & 2 & 0 \\
\hline \hline
\mbox{Facet} & p_1 & q_1 & p_2 & q_2 & p_3 & q_3 & p_4 & q_4 & & \\
\hline
\mbox{index} & 1 & 0 & 0 & 1 & 0 & 0 & 1 & 1 & & \\
\hline
\end{array}$$

  The total number of wedges done equals $23$, hence we get two polytopes of 
dimension $47$ with $51$ facets. Their moment-angle manifolds have dimension 
$98$.

  We can compute their Betti numbers, what characterizes them, up to 
diffeomorphism, and we assert they are the same in both cases. We indeed find 
that both moment-angle manifolds are diffeomorphic to the following connected 
sum of sphere products:
$${\# \atop 2} (S^{37} \times S^{61}) \# (S^{38} \times S^{60})
{\# \atop 17} (S^{39} \times S^{59}) {\# \atop 19} (S^{40} \times S^{58})
{\# \atop 20} (S^{41} \times S^{57}) {\# \atop 22} (S^{42} \times S^{56})
{\# \atop 38} (S^{43} \times S^{55})$$
$${\# \atop 27} (S^{44} \times S^{54}){\# \atop 29} (S^{45} \times S^{53}) 
{\# \atop 51} (S^{46} \times S^{52}){\# \atop 54} (S^{47} \times S^{51}) 
{\# \atop 75} (S^{48} \times S^{50}){\# \atop 36} (S^{49} \times S^{49})$$

  We also assert that these polytopes have different Betti numbers. We know 
that $b^{-1,0} = b^{46,51} = 1$. All other nonzero ones appear in the following 
table (when two numbers are given, the first one is the one of $P_1 $, the 
second the one of $P_2 $):

$$\begin{array}{|c||c|c|c|c|c|c|}
\hline
\mbox{Bet. num.} & b^{17,19} , b^{28,32} & b^{18,20} , b^{27,31} & 
b^{19,21} , b^{26,30} & b^{20,22} , b^{25,29} & b^{21,23} , b^{24,28} & 
b^{22,24} , b^{23,27} \\
\hline
\mbox{value} & 2 & 17 & 15 & 30,29 & 22,25 & 40,37 \\
\hline \hline
\mbox{Bet. num.} & b^{17,20} , b^{28,31} & b^{18,21} , b^{27,30} & 
b^{19,22} , b^{26,29} & b^{20,23} , b^{25,28} & b^{21,24} , b^{24,27} & 
b^{22,25} , b^{23,26} \\
\hline
\mbox{value} & 1 & 19 & 22 & 27 & 51 & 75 \\
\hline \hline
\mbox{Bet. num.} & b^{23,25} , b^{22,26} & b^{24,26} , b^{21,25} & 
b^{25,27} , b^{20,24} & b^{26,28} , b^{19,23} & b^{27,29} , b^{18,22} & 
b^{28,30} , b^{17,21} \\
\hline
\mbox{value} & 36 & 14,17 & 7,4 & 8,9 & 5 & 0 \\
\hline
\end{array}$$

  We immediately check this is coherent with the diffeomorphism type of their 
moment-angle manifold.

  We moreover can count the number of vertices of each polytope. The polytope 
$P_1 $ has two vertices more ($33686$) than $P_2 $ ($33684$). This proves, 
as announced, that these two polytopes have different $f$-vectors.

\vskip 10mm

{\footnotesize {Bosio Fr\'ed\'eric \\
Universit\'e de Poitiers \\
UFR Sciences SP2MI \\
D\'epartement de Math\'ematiques \\
UMR CNRS 7348 \\
Teleport 2 \\
Boulevard Marie et Pierre Curie \\
BP 30179 \\
86962 Futuroscope Chasseneuil CEDEX 

e-mail~: Frederic.Bosio@math.univ-poitiers.fr}}

\end{document}